\documentclass[12pt,twoside]{article} 
\usepackage{amsmath,amssymb,mathrsfs,graphicx}

\newtheorem{theorem}{Theorem}[section]
\newtheorem{lemma}[theorem]{Lemma}
\newtheorem{proposition}[theorem]{Proposition}

\newtheorem{definition}[theorem]{Definition}

\newtheorem{rmrk}[theorem]{Remark}

\DeclareMathAlphabet{\mathbfit}{OML}{cmm}{b}{it}

\makeatletter
\@addtoreset{equation}{section}
\makeatother

\newenvironment{remark}
{\begin{rmrk} \em}
{\end{rmrk}}

\newcommand{\fn} {function}
\newcommand{\me} {measure}

\newcommand{\erg} {ergodic}
\newcommand{\sy} {system}

\newcommand{\pr} {probability}
\newcommand{\dsy} {dynamical system}
\renewcommand{\o} {orbit}

\newcommand{\R} {\mathbb{R}}

\newcommand{\Z} {\mathbb{Z}}
\newcommand{\N} {\mathbb{N}}
\newcommand{\qed} {\hfill {\small Q.E.D.} \par\medskip}
\newcommand{\skippar} {\par\medskip}
\newcommand{\ds} {\displaystyle}
\newcommand{\proof} {\noindent \textsc{Proof.} }
\newcommand{\proofof}[1] {\noindent \textsc{Proof of {#1}.} }
\newcommand{\article}[3] {\textsc{{#1}}, {\itshape {#2}}, {{#3}}.}
\newcommand{\book}[3] {\textsc{{#1}}, {\itshape {#2}}, {{#3}}.}
\newcommand{\vol} {\textbf}
\newcommand{\eps} {\varepsilon}

\newcommand{\rset}[2] {\left\{ #1 \: \left| \: #2 \right. \! \right\} }

\newcommand{\symmdiff} {\triangle}
\newcommand{\la} {\langle}
\newcommand{\ra} {\rangle}
\newcommand{\into} {\longrightarrow}

\renewcommand{\iff} {if and only if\ }

\newcommand{\m} {mixing}
\newcommand{\ob} {observable}
\newcommand{\ps} {X}            
\newcommand{\sca} {\mathscr{A}}           
\newcommand{\scb} {\mathscr{B}}           
\newcommand{\scn} {\mathscr{N}}           
\newcommand{\scv} {\mathscr{V}}
\newcommand{\go} {\mathcal{G}}            
\newcommand{\lo} {\mathcal{L}}            
\newcommand{\iv} {V \nearrow \ps}
\newcommand{\ivlim} {\lim_{\iv}}
\newcommand{\ivlimtwo} {\lim_{{\iv} \atop {n \to \infty}}} 
\newcommand{\avg} {\overline{\mu}}
\newcommand{\jv} {\mathbb{J}_V}           
\newcommand{\eo} {\mathcal{E}}            
\newcommand{\E} {\mathbb{E}}

\begin{document}

\title{\textbf{Exactness, K-property \\ and infinite mixing}}

\author{\textsc{Marco Lenci}
\thanks{
Dipartimento di Matematica, 
Universit\`a di Bologna, 
Piazza di Porta S.\ Donato 5, 
40126 Bologna, Italy.
E-mail: \texttt{marco.lenci@unibo.it} } 
\thanks{Istituto Nazionale di Fisica Nucleare, 
Sezione di Bologna,
Via Irnerio 46,
40126 Bologna, Italy.}
}

\date{\normalsize
Version published in \\ 
\emph{Publicaciones Matem\'aticas del Uruguay} 
\textbf{14} (2013), 159-170 \\
+ minor fixes \\
\vspace{8pt}
November 2013
}

\maketitle

\begin{abstract}
  We explore the consequences of exactness or K-mixing on the notions
  of mixing (a.k.a.\ \emph{infinite-volume mixing}) recently devised
  by the author for infinite-measure-preserving dynamical systems.

\end{abstract}

\section{Introduction}
\label{sec-intro}

Currently, in infinite \erg\ theory, there is a renewed interest in
the issues related to \m\ for infinite-\me-preserving (or just
nonsingular) \dsy s, in short \emph{infinite \m} (see \cite{z, ds,
  limix, i3, dr, mt, lp, a2, ko, t1}, and some applications in
\cite{i1, i2, amps, liutam, t2}).

The present author recently introduced some new notions of infinite
\m, based on the concept of \emph{global \ob} and
\emph{infinite-volume average} \cite{limix}. In essence, a global \ob\
for an infinite, $\sigma$-finite, \me\ space $(\ps, \sca, \mu)$ is
\fn\ in $L^\infty (\ps, \sca, \mu)$ that ``looks qualitatively the
same'' all over $\ps$. This is in contrast with a \emph{local \ob}, 
whose support is essentially localized, so that the \fn\ is integrable.

Postponing the mathematical details to Section \ref{sec-setup}, the
purpose of the global \ob s is basically twofold.  First,
the past attempts to a general definition of infinite mixing involved
mainly local \ob s (equivalently, finite-\me\ sets), and the problems
with such definitions seemed to depend on that.  Second, seeking
inspiration in statistical mechanics (which is the discipline of
mathematical physics that has successfully dealt with the question of
predicting measurements in very large, formally infinite, \sy s), one
realizes that many quantities of interest are \emph{extensive \ob s},
that is, objects that behave qualitatively in the same way in
different regions of the phase space. (More detailed discussions about
these points are found in in \cite{limix, liutam}.)

Extensive \ob s are ``\me d'' by taking averages over large portions
of the phase space. We import that concept too, by defining the
infinite-volume average of a global \ob\ $F: \ps \into \R$ as
\begin{equation}
  \label{avg-intro}
  \avg(F) := \ivlim\, \frac1 {\mu(V)} \int_V F \, d\mu.
\end{equation}
Here $V$ is taken from a family of ever larger but finite-\me\ sets
that somehow covers, or \emph{exhausts} the whole of $\ps$. The
precise meaning of the limit above will be given in Section
\ref{sec-setup}.

Now, let us consider a \me-preserving \dsy\ on $(\ps, \sca, \mu)$.
For the sake of simplicity, let us restrict to the discrete-time case:
this means that we have a measurable map $T: \ps \into \ps$ that
preserves $\mu$. Choosing two suitable classes of global and local \ob
s, respectively denoted $\go$ and $\lo$, we give five definitions of
infinite mixing. These fall in two categories, exemplified as follows.

Using the customary (abuse of) notation $\mu(g) := \int_\ps g\,
d\mu$. we say that the \sy\ exhibits:
\begin{itemize}
\item \emph{global-local mixing} if, $\forall F \in \go$, $\forall g
  \in \lo$, $\ds \lim_{n \to \infty} \, \mu((F \circ T^n) g) = \avg(F)
  \mu(g)$;
  
\item \emph{global-global mixing} if, $\forall F,G \in \go$, $\ds
  \lim_{n \to \infty} \, \avg ( (F \circ T^n) G ) = \avg(F) \,
  \avg(G)$.
\end{itemize}

Disregarding for the moment the mathematical issues connected to the
above notions, we focus on the interpretation of global-local
mixing. Restricting, without loss of generality, to local \ob s $g \ge
0$ with $\mu(g) = 1$, and defining $d\mu_g := g \, d\mu$, the above
limit reads:
\begin{equation}
  \label{push-forward}
  \lim_{n \to \infty} T_*^n \mu_g (F) = \avg(F),
\end{equation}
where the \me\ $T_*^n \mu_g$ is the push-forward of $\mu_g$ via the
dynamics $T^n$ (in other words, $T_*^n \mu_g := \mu_g \circ T^{-n} =
\mu_{P^n g}$, where $P$ is the Perron-Frobenius operator relative to
$\mu$, cf.\ (\ref{coupling})-(\ref{def-pf})). If (\ref{push-forward})
occurs for all $g \in \lo$ and $F \in \go$, the above is a sort of
``convergence to equilibrium'' for all \emph{initial states} given by
$\mu$-absolutely continuous \pr\ \me s. In this sense the functional
$\avg$ (not a \me!)\ plays the role of the \emph{equilibrium state}.

Exactness and K-mixing (a.k.a.\ the K-property) are notions that exist
and have the same definition both in finite and infinite \erg\
theory. In finite \erg\ theory they are known to be very strong
properties, as they imply mixing of all orders, cf.\ definition
(\ref{mix-all-o}).  The purpose of this note is to explore their
implications in terms of the notions of infinite \m\ introduced in
\cite{limix}.

As we will see below (Theorem \ref{thm-main}\emph{(a)}), the most
notable of such implications is a weak form of global-local mixing,
whereby any pair of \me s $\mu_g$, $\mu_h$, as introduced earlier, 
are \emph{asymptotically coalescing}, in the sense that
\begin{equation}
  \lim_{n \to \infty} \left( T_*^n \mu_g (F) - T_*^n \mu_h (F) \right) = 0,
\end{equation}
for all $F \in \go$.

In the next section we review the five definitions of global-local and
global-global \m, together with the already known (though with a
different name) definition of local-local mixing. In Section
\ref{sec-ek} we prepare, state and prove Theorem \ref{thm-main}, which
lists some consequences of exactness and the K-property. Finally, in
Section \ref{sec-eo}, we introduce the space of the \emph{equilibrium
  \ob s}, which is a purely \erg-theoretical construct in which some
information about global-local \m\ can be recast.

\section{Definitions of infinite mixing}
\label{sec-setup}

Let $(\ps, \sca, \mu, T)$ be a \me-preserving \dsy, where $(\ps,
\sca)$ is a \me\ space, $\mu$ an infinite, $\sigma$-finite, \me\ on
it, and $T$ a $\mu$-endomorphism, that is, a measurable surjective map
that preserves $\mu$ (i.e., $\mu(T^{-1} A) = \mu(A)$, $\forall A \in
\sca$).

Denoting by $\sca_f := \rset{A \in \sca} {\mu(A) < \infty}$ the class
of finite-\me\ sets, we assume that the following additional structure
is given for the \dsy:

\begin{itemize}
\item A class of sets $\scv \subset \sca_f$, called the {\bfseries
  exhaustive family}. The elements of $\scv$ will be generally
  indicated with the letter $V$.

\item A subspace $\go \subset L^\infty (\ps, \sca, \mu; \R)$, whose
  elements are called the {\bfseries global \ob s}. These \fn s are
  indicated with uppercase Roman letters ($F,G$, etc.).

\item A subspace $\lo \subset L^1(\ps, \sca, \mu; \R)$ whose elements
  are called the {\bfseries local \ob s}. These \fn s will be
  indicated with lowercase Roman letters ($f,g$, etc.).
\end{itemize}
A discussion on the role and the choice of $\scv, \go, \lo$ is given
in $\cite{limix}$, together with the proofs of most assertions made in
this section.

We assume that $\scv$ contains at least one sequence $(V_j)_{j \in
\N}$, ordered by inclusion, such that $\bigcup_j V_j = X$. (In
actuality, this requirement is never used in the proofs, but, since
the elements of $\scv$ are regarded as large and ``representative''
regions of the phase space $\ps$, we keep it to give ``physical''
meaning to the concept of infinite-volume average, see below.)  We
also assume that $1 \in \go$ (with the obvious notation $1(x) := 1$,
$\forall x \in \ps$).

\begin{definition}
  \label{def-iv}
  Let $\scv$ be the aforementioned exhaustive family. For $\phi: \scv
  \into \R$, we write
  \begin{displaymath}
    \ivlim \phi(V) = \ell
  \end{displaymath}
  when  
  \begin{displaymath}
    \lim_{M \to \infty} \: \sup_{{V \in \scv} \atop {\mu(V) \ge M}} | 
    \phi(V) - \ell | = 0.
  \end{displaymath}
  We call this the `$\mu$-uniform infinite-volume limit w.r.t.\ the
  family $\scv$', or, for short, the {\bfseries infinite-volume
  limit}.
\end{definition}

We assume that, $\forall n \in \N$,
\begin{equation}
 \label{a-vol}
 \mu( T^{-n} V \symmdiff V ) = o(\mu(V)), \quad \mbox{as } \iv.
\end{equation}
This is reasonable because, if a large $V \in \scv$ is to be
considered a finite-\me\ substitute for $\ps$, it makes sense to
require that a finite-time application of the dynamics does not change
it much. Finally, the most crucial assumption is that,
\begin{equation}
  \label{a-avg}
  \forall F \in \go, \quad \exists \, \avg(F) := \ivlim\, \frac1 
  {\mu(V)} \int_V F \, d\mu.
\end{equation}
$\avg(F)$ is called the {\bfseries infinite-volume average} of $F$
w.r.t.~$\mu$. It easy to check that $\avg$ is $T$-invariant, i.e., for
all $F \in \go$ and $n \in \N$, $\avg(F \circ T^n)$ exists and equals
$\avg(F)$ \cite{limix}.

\skippar

With this machinery, we can give a number of definitions of infinite
\m\ for the \dsy\ $(\ps, \sca, \mu, T)$ endowed with the
\emph{structure of \ob s} $(\scv, \go, \lo)$.

The following three definitions will be called {\bfseries global-local
\m}, as they involve the coupling of a global and a local \ob.  We
say that the \sy\ is \m\ of type
\begin{description}
\item[(GLM1)] \, if, $\forall F \in \go$, $\forall g \in \lo$ with
  $\ds \mu(g)=0$, $\ds \lim_{n \to \infty} \, \mu((F \circ T^n) g) =
  0$;
  
\item[(GLM2)] \, if, $\forall F \in \go$, $\forall g \in \lo$, $\ds
  \lim_{n \to \infty} \, \mu((F \circ T^n) g) = \avg(F) \mu(g)$;

\item[(GLM3)] \, if, $\forall F \in \go$, $\ds \lim_{n \to \infty} \,
  \sup_{g \in \lo \setminus 0 } \| g \|_1^{-1} \left| \mu((F \circ
  T^n) g) - \avg(F) \mu(g) \right| = 0$,
\end{description}
where $\| \cdot \|_1$ is the norm of $L^1(\ps, \sca, \mu; \R)$.

Clearly, \textbf{(GLM1--3)} are listed in increasing order of
strength, with \textbf{(GLM2)} being possibly the most natural
definition one can give for the time-decorrelation between a global
and a local \ob\ (recall that $\avg(F \circ T^n) = \avg(F)$).
\textbf{(GLM3)} is a uniform version of it, with important
implications (cf.\ Proposition \ref{prop-mgl3-mgg2}), while
\textbf{(GLM1)} is a much weaker version, as will become apparent in
the remainder.

\skippar

Although this note is mostly concerned with global-local \m, one can
also consider the decorrelation of two global \ob s, namely {\bfseries
global-global \m}. For this we need the following terminology:

\begin{definition}
  \label{def-ivtwo}
  For $\scv$ as defined above and $\phi: \scv \times 
  \N \into \R$, we write 
  \begin{displaymath}
    \ivlimtwo \phi(V,n) = \ell
  \end{displaymath}
  to mean
  \begin{displaymath}
    \lim_{M \to \infty} \sup_{{V \in \scv} \atop 
    {{\mu(V) \ge M} \atop {n \ge M}}} | \phi(V,n) - \ell | = 0.
  \end{displaymath}
  As $n$ will take the role of time, we refer to this limit as the
  `joint infinite-volume and time limit'.
\end{definition}

For $F \in L^\infty$ and $V \in \scv$, let us also denote $\mu_V (F)
:= \mu(V)^{-1} \int_V F d\mu$.  We say that the \sy\ is \m\ of type
\begin{description}
\item[(GGM1)] \, if, $\forall F,G \in \go$, $\ds \lim_{n \to \infty}
  \, \avg ( (F \circ T^n) G ) = \avg(F) \, \avg(G)$;

\item[(GGM2)] \, if, $\forall F,G \in \go$, $\ds \ivlimtwo \mu_V ((F
  \circ T^n) G) = \avg(F) \, \avg(G)$.
\end{description}

Though \textbf{(GGM1)} seems the cleaner of the two versions, it has
the serious drawback that, for $n \in \N$, $\avg ( (F \circ T^n) G )$
might not even exist, for there is no provision in our hypotheses to
guarantee the ring property for condition (\ref{a-avg}) (namely,
$\exists \avg(F), \avg(G) \Rightarrow \exists \avg(FG))$. Nor do we
want one, if we are to keep our framework general enough.
\textbf{(GGM2)} solves this question of wellposedness, and is in some
sense stronger than \textbf{(GGM1)}:

\begin{proposition}
  \label{prop-mgg2-1}
  If $F,G \in \go$ are such that $\avg ( (F \circ T^n) G )$ exists for
  all $n$ large enough (depending on $F,G$), then
  \begin{equation}
    \label{mgg2-1-impl}
    \ivlimtwo \mu_V ((F \circ T^n) G) = \ell \quad \Longrightarrow 
    \quad \lim_{n \to \infty} \avg ((F \circ T^n) G) = \ell.
  \end{equation}
  In particular, if the above hypothesis holds $\forall F,G \in \go$,
  then {\rm\bfseries (GGM2)} implies {\rm\bfseries (GGM1)}.
\end{proposition}

\proof From Definition \ref{def-ivtwo}, the left limit of
(\ref{mgg2-1-impl}) implies that, $\forall \eps > 0$, $\exists M =
M(\eps)$ such that
\begin{equation}
  \label{mgg2-1-10}
  \ell - \eps \le \mu_V ((F \circ T^n) G) \le \ell + \eps 
\end{equation}
for all $V \in \scv$ with $\mu(V) \ge M$ and all $n \ge M$. By
hypothesis, if $M$ is large enough, the infinite-volume limit of the
above middle term exists $\forall n \ge M$ and equals $\avg ( (F \circ
T^n) G )$. Upon taking such limit, what is left of (\ref{mgg2-1-10})
and its conditions of validity is the very definition of the right
limit in (\ref{mgg2-1-impl}).
\qed

With reasonable hypotheses on the structure of $\go$ and $\lo$, the
strongest version of global-local \m\ implies the ``strongest''
version of global-global \m.  The following proposition is a
simplified version of a similar result of \cite{limix} (for an
intuitive understanding of the hypotheses, see Proposition 3.2 and
Remark 3.3 there).

\begin{proposition}
  \label{prop-mgl3-mgg2}
  Suppose there exist a family $( \psi_j )_{j \in \N}$ of real-valued
  \fn s of $\ps$ (this will play the role of a partition of unity) and
  a family $( \jv )_{V \in \scv}$ of finite subsets of $\N$ such that:
  \begin{enumerate}
  \item[(i)] $\forall j \in \N$, $\psi_j \ge 0$;
  \item[(ii)] $\forall G \in \go$, $\forall j \in \N$, $G\psi_j \in
    \lo$;
  \item[(iii)] in the limit $\iv$, $\ds \left\| \sum_{j \in \jv}
    \psi_j - 1_V \right\|_1 = o(\mu(V))$,
  \end{enumerate}
  where $1_V$ is the indicator \fn\ of $V$. Then {\rm\bfseries
  (GLM3)} implies {\rm\bfseries (GGM2)}.
\end{proposition}

\proofof{Proposition \ref{prop-mgl3-mgg2}} Since the limit in
{\bfseries (GGM2)} is trivial when $G$ is a constant, and since the
global \ob s are bounded \fn s, it is no loss of generality to prove
{\bfseries (GGM2)} for the case $G \ge 0$ only.

The proof follows upon verification that the \fn s $g_j := G\psi_j$
verify all the hypotheses of Proposition 3.2 of \cite{limix} (cf.\
also Remark 3.3). Notice that the identity $G = \sum_j g_j$ (which
makes sense insofar as $( \psi_j )_j$ is a partition of unity) is
illustrative and not really used in the proof there.  \qed

\skippar 

Since the five definitions presented above deal with the decorrelation
of, first, a global and a local \ob, and then two global \ob s,
symmetry considerations would induce one to give a definition of
{\bfseries local-local \m} as well. A reasonable possibility would be
to call a \dsy\ \m\ of type
\begin{description}
\item[(LLM)] \, if, $\forall f \in \lo \cap \go$, $g \in \lo$,
  $\ds \lim_{n \to \infty} \, \mu ( (f \circ T^n) g ) = 0$.
\end{description}
In fact, this definition already exists, as it is easy to check that,
in the most general case (that is, $\go = L^\infty$, $\lo = L^1$), a
\dsy\ is {\bfseries (LLM)} if and only if, $\forall A,B \in \sca_f$,
$\lim_{n \to \infty} \mu(T^{-n} A \cap B) = 0$, i.e., if and only if
the \sy\ is of \emph{zero type} \cite{hk} (cf.\ also \cite{ds,
ko}). Incidentally, this is the same definition that Krengel and
Sucheston call `\m', for an infinite-\me-preserving \dsy\ \cite{ks}.

\section{Exactness and K-property}
\label{sec-ek}

Two of the few definitions that are copied verbatim from finite to
infinite \erg\ theory are those of exactness and K-\m. Though they are
well known, we repeat them here for completeness. We state the
versions for \me-preserving maps, but they can be given for
nonsingular maps as well ($T$ is nonsingular if $\mu(A) =0 \Rightarrow
\mu(T^{-1} A) =0$).

Let us denote by $\scn$ the \emph{null $\sigma$-algebra}, i.e., the
$\sigma$-algebra that only contains the zero-\me\ sets and their
complements. Also, given two $\sigma$-algebras $\sca, \scb$, we write
$\sca = \scb$ mod $\mu$ if $\forall A \in \sca$, $\exists B \in \scb$
with $\mu( A \symmdiff B ) = 0$, and viceversa; equivalently, the
$\mu$-completions of $\sca$ and $\scb$ are the same.

\begin{definition}
  \label{def-exact}
  The \me-preserving \dsy\ $(\ps, \sca, \mu, T)$ is called {\bfseries
  exact} if
  \begin{displaymath}
    \bigcap_{n=0}^\infty T^{-n} \sca = \scn \ \mathrm{mod} \ \mu.
  \end{displaymath}
\end{definition}

Since exactness implies that $T^{-1} \sca \ne \sca$ mod $\mu$, a
nontrivial exact $T$ cannot be an automorphism of the \me\ space
$(\ps, \sca, \mu)$---although in some sense an invertible map can
still be exact, cf.\ Remark \ref{rk-inv-exact} below.

The counterpart of exactness for automorphisms is the following:

\begin{definition}
  \label{def-kmix}
  The invertible \me-preserving \dsy\ $(\ps, \sca, \mu, T)$ possesses
  the {\bfseries K-property} (from A.~N.~Kolmogorov) if $\exists \scb
  \subset \sca$ such that:
  \begin{itemize}
  \item[(i)] $\scb \subset T \scb$;
  \item[(ii)] $\ds \bigvee_{n=0}^\infty T^n \scb = \sca \ \mathrm{mod}
    \ \mu$;
  \item[(iii)] $\ds \bigcap_{n=0}^\infty T^{-n} \scb = \scn \
    \mathrm{mod} \ \mu$.
  \end{itemize} 
  In this case, one also says that the \dsy\ is K-\m, or that $T$ is a
  K-automorphism of $(\ps, \sca, \mu)$.
\end{definition}

\begin{remark}
  \label{rk-inv-exact}
  Comparing Definition \ref{def-exact} with condition \emph{(iii)} of
  Definition \ref{def-kmix}, one might be tempted to say that, if
  $(\ps, \sca, \mu, T)$ has the K-property, then $(\ps, \scb, \mu, T)$
  is exact. This is not \emph{technically} correct because, in all
  nontrivial cases, the inclusion in Definition
  \ref{def-kmix}\emph{(i)} is strict, thus $T$ is not a self-map of
  the \me\ space $(\ps, \scb, \mu)$.  That said, if $(\ps, \sca, \mu)$
  is a Lebesgue space, $(\ps, \scb, \mu, T)$ is still \emph{morally}
  exact, in the following sense. Assume w.l.g.\ that $\scb$ is
  complete, let $\ps_\scb$ be the measurable partition of $\ps$ that
  generates $\scb$. (In a Lebesgue space there is a one-to-one
  correspondence, modulo null sets, between complete
  sub-$\sigma$-algebras and measurable partitions \cite{r}.)  $\scb$
  can be lifted to a $\sigma$-algebra for $\ps_\scb$, which we keep
  calling $\scb$. Also, defining $T_\scb ([x]) := [T(x)]$ (where $[x]$
  denotes the element of $\ps_\scb$ that contains $x$), we verify that
  $T_\scb$ is well defined as a self-map of $(\ps_\scb, \scb, \mu)$
  (in fact, from Definition \ref{def-kmix}\emph{(i)}, $\ps_{T \scb}$
  is a sub-partition of $\ps_\scb$) and $T_\scb^{-1} A = T^{-1} A$,
  $\forall A \in \scb$ (with the understandable abuse of notation
  whereby $A$ denotes both a subset of $\ps_\scb$ and a subset of
  $\ps$). This and Definition \ref{def-kmix}\emph{(iii)} show that
  $T_\scb$ is an exact endomorphism of $(\ps_\scb, \scb, \mu)$.  Of
  course, in all of the above, $\scb$ can be replaced by $\scb_m :=
  T^m \scb$, for all $m \in \Z$ (because $\scb_m$ can be used in lieu
  of $\scb$ in Definition \ref{def-kmix}).
\end{remark}

In finite \erg\ theory, both exactness and the K-property imply
\emph{\m\ of all orders}, namely, $\forall k \in \Z^+$ and $A_1, A_2,
\ldots, A_k \in \sca$,
\begin{equation}
  \label{mix-all-o}
  \mu(A_1 \cap T^{-n_2} A_2 \cap \cdots T^{-n_k} A_k) \into
  \mu(A_1) \mu(A_2) \cdots \mu(A_k),
\end{equation}
whenever $n_2 \to \infty$ and $n_{i+1} - n_1 \to \infty$, $\forall
i=2, \ldots k-1$ \cite{q}. (In (\ref{mix-all-o}) we have assumed
$\mu(\ps) = 1$.)

One would expect such strong properties to have consequences also in
infinite \erg\ theory. This is the case, as we describe momentarily.
But first we need some elementary formalism from the functional
analysis of \dsy s.  For $F \in L^\infty$ and $g \in L^1$, let us
denote
\begin{equation}
  \label{coupling}
  \la F,g \ra := \mu(Fg).
\end{equation}
Define the \emph{Koopman operator} $U: L^\infty \into L^\infty$ as $UF
:= F \circ T$. Its adjoint for the above coupling is called the
\emph{Perron-Frobenius operator}, denoted $P: L^1 \into L^1$. Its
defining identity is
\begin{equation}
  \label{def-pf}
  \la UF, g \ra = \la F, Pg \ra.
\end{equation}
Let us explain in detail how $P$ is defined through (\ref{def-pf}).
Take $g \in L^1$ and assume for the moment $g \ge 0$. Take also $F =
1_A$, with $A \in \sca$. We see that $\la UF, g \ra = \int_{T^{-1} A}
g \, d\mu$. Since $T$ preserves $\mu$ and is thus nonsingular w.r.t.\
it, and since the \me\ space is $\sigma$-finite, the Radon-Nykodim
Theorem yields a locally-$L^1$, positive, \fn\ $Pg: \ps \into \R$ such
that $ \int_{T^{-1} A} g \, d\mu = \int_A (Pg) d\mu = \la F, Pg
\ra$. Using $F = 1_\ps = 1$, we see that $Pg \in L^1$ with $\| Pg \|_1
= \| g \|_1$. For a general $g \in L^1$, we write $g = g^+ - g^-$,
where $g^+$ and $g^-$ are, respectively, the positive and negative
parts of $g$. Then $Pg := Pg^+ - Pg^-$ is also in $L^1$ and
\begin{equation}
  \| Pg \|_1 \le \| g \|_1.
\end{equation}
Therefore, through approximations of $F$ via simple \fn s (in the 
$L^\infty$-norm), one can extend (\ref{def-pf}) to all $F \in L^\infty$.

In the process, we have learned that $P$ is a positive operator ($g
\ge 0 \Rightarrow Pg \ge 0$) and $\| P \| =1$, whereas, obviously, $U$
is a positive isometry.  Moreover, it is easy to see that $Pg = g$,
with $g \ge 0$, \iff $g$ is an invariant density, i.e., if $\mu_g$
defined by $d \mu_g / d\mu = g$ is an invariant \me. (In fact, had we
defined (\ref{coupling}) for $F \in L^1$ and $g \in L^\infty$,
(\ref{def-pf}) would have defined a positive operator $P: L^\infty
\into L^\infty$, with $\| P \| =1$, and such that $P1 = 1$.)

\skippar

Most of the remainder of this note will be based on an important
theorem by Lin \cite{li} (see also \cite{a} for a nice short proof).

\begin{theorem}
  \label{thm-lin}
  The nonsingular \dsy\ $(\ps, \sca, \mu, T)$ is exact if and only if,
  $\forall g \in L^1$ with $\mu(g) = 0$, $\ds \lim_{n \to \infty} \|
  P^n g \|_1 = 0$.
\end{theorem}

In the rest of the paper we assume to be in one of the following two
cases:

\begin{itemize}
\item[(H1)] $(\ps, \sca, \mu, T)$ is exact.  $\scv$ is any exhaustive
  family that verifies (\ref{a-vol}). $\go = L^\infty$. $\lo = L^1$.
  (Given the assumptions of Section \ref{sec-setup}, this corresponds
  to the most general choice of $\scv, \go, \lo$.)
  
\item[(H2)] $(\ps, \sca, \mu, T)$ is K-\m\ (thus $T$ is an
  automorphism). $\scv$ is any exhaustive family that verifies
  (\ref{a-vol}). $\go$ is the closure, in $L^\infty$, of
  $\bigcup_{m>0} L^\infty(\scb_m)$, where $\scb_m = T^m \scb$, as
  defined in Remark \ref{rk-inv-exact}. Lastly, $\lo = L^1$.
\end{itemize}

\begin{theorem}
  \label{thm-main}
  Under either {\rm (H1)} or {\rm (H2)},
  \begin{itemize}
  \item[(a)] {\rm\bfseries (GLM1)} holds true;
  \item[(b)] {\rm\bfseries (LLM)} holds true;   
  \item[(c)] {\rm\bfseries (GGM2)} implies {\rm\bfseries (GLM2)};
  \item[(d)] If, $\forall F \in \go$, $\exists g_F \in \lo$, with 
    $\mu(g_F) \ne 0$, such that
    \begin{displaymath}
      \lim_{n \to \infty} \, \mu((F \circ T^n) g_F) = \avg(F) \mu(g_F),
    \end{displaymath}
    then {\rm\bfseries {(GLM2)}} holds true.
  \end{itemize}
\end{theorem}

As anticipated in the introduction, {\bfseries (GLM1)} (which is the
most important assertion of the theorem) means that the evolutions of
two absolutely continuous initial \me s become indistinguishable, as
time goes to infinity. We may call this phenomenon \emph{asymptotic
coalescence}. This implies that they will return the same
measurements of global \ob s, but not that this measurements will
converge (in which case we would have a sort of convergence to
equilibrium). In fact, for many interesting \sy s, it is not hard to
construct $F \in L^\infty$ such that $\la F, P^n g \ra$ does not
converge for all $g \in L^1$.

This is not surprising, for, even in finite \erg\ theory, certain
proofs of mixing, or decay of correlation, are divided in two parts:
asymptotic coalescence and the convergence of \emph{one} initial
\me. The difference there is that the latter is usually easy.

The remainder of this section is devoted to the following:

\skippar

\proofof{Theorem \ref{thm-main}} Let us start by proving assertion
\emph{(a)}, namely {\bfseries (GLM1)}. We use the formalism of
functional analysis outlined earlier in the section.

If (H1) is the case, the proof is immediate: for $F \in L^\infty$ and
$g \in L^1$, with $\mu(g) = 0$,
\begin{equation}
  \label{main-10}
  | \mu((F \circ T^n) g) | = | \la F , P^n g \ra | \le \| F \|_\infty \, 
  \| P^n g \|_1 \to 0,
\end{equation}
as $n \to \infty$, by Theorem \ref{thm-lin}.

In the case (H2), let us observe that, by easy density arguments, all
the definitions {\bfseries (GLM1--3)} hold true if they are verified
w.r.t.\ $\go'$ and $\lo'$ which are subspaces of $\go$ and $\lo$,
respectively, in the $L^\infty$- and $L^1$-norms. We can take $\go' :=
\bigcup_{m>0} L^\infty(\scb_m)$ (which is dense in $\go$ by
definition) and $\lo' := \bigcup_{m>0} L^1(\scb_m)$, which is dense in
$\lo = L^1(\sca)$ by the K-property \cite{a}.  Therefore, it suffices
to show {\bfseries (GLM1)} for a general $m > 0$ and $\forall F \in
L^\infty (\scb_m)$, $\forall g \in L^1 (\scb_m)$ with $\mu(g) = 0$.

Using the arguments and the notation of Remark \ref{rk-inv-exact}, we
denote by $\hat{F}$ the \fn\ induced by $F$ on $\ps_{\scb_m}$ (i.e.,
$\hat{F}([x]) := F(x)$), and analogously for all the other
$\scb_m$-measurable \fn s. We observe that $F \circ T^n$ is
$\scb_m$-measurable and $\widehat{F \circ T^n} = \hat{F} \circ
T_{\scb_m}^n$. Thus
\begin{equation}
  \label{main-15}
  \mu((F \circ T^n) g)  = \mu( (\hat{F} \circ T_{\scb_m}^n) 
  \hat{g} ),
\end{equation}
where the r.h.s.\ is regarded as an integral in $\ps_{\scb_m}$.  Since
$(\ps_{\scb_m}, \scb_m, \mu, T_{\scb_m})$ is exact, and $\mu(\hat{g})
= \mu(g) = 0$, we use (\ref{main-15}) in (\ref{main-10}) to prove that
the l.h.s.\ of (\ref{main-15}) vanishes, as $n \to \infty$.

The following is a corollary of {\bfseries (GLM1)}.

\begin{lemma}
  \label{lem-main}
  Assume either {\rm (H1)} or {\rm (H2)}, and fix $F \in \go$.  If,
  for some $\ell \in \R$ and $\eps \ge 0$, the limit
  \begin{displaymath}
    \limsup_{n \to \infty} \left| \frac{ \mu((F \circ T^n) g) } {\mu(g)} 
    - \ell \right| \le \eps
  \end{displaymath}
  holds for some $g \in \lo$ (with $\mu(g) \ne 0$), then it holds for
  all $g \in \lo$ (with $\mu(g) \ne 0$).
\end{lemma}

\proofof{Lemma \ref{lem-main}} Suppose the above limit holds for $g_0
\in \lo$. Take any other $g \in \lo$, with $\mu(g) \ne 0$. We have:
\begin{equation}
\begin{split}
  & \left| \frac{ \mu((F \circ T^n) g) } { \mu(g) } - \ell \right| \\
  \le & \left| \mu \! \left( (F \circ T^n) \left( \frac{g}{\mu(g)} -
  \frac{g_0}{\mu(g_0)} \right) \right) \right| + \left| 
  \frac{ \mu((F \circ T^n) g_0) } { \mu(g_0) } - \ell \right| .
\end{split}
\end{equation}
By {\bfseries (GLM1)}, the first term of the above r.h.s.\ vanishes 
as $n \to \infty$, whence the assertion.
\qed

Going back to the proof of Theorem \ref{thm-main}, we see that Lemma
\ref{lem-main} immediately implies assertion \emph{(d)}.

As for \emph{(b)}, again we prove it for both cases (H1) and (H2) at
the same time.  W.l.g., let us assume that $\sca \ne \scn$ mod $\mu$
(otherwise $L^1$ would be trivial). We claim that
\begin{equation}
  \label{main-20}
  \sup_{A \in \sca_f} \mu(A) = \infty.
\end{equation}
In fact, since $\sca$ is not trivial, the above sup is positive. If it
equalled $M \in \R^+$, it would be easy to construct an invariant set
$B$ with $0 < \mu(B) \le M$. But $\mu(\ps) = \infty$, therefore $T$
would not be \erg, contradicting both (H1) and (H2).

Now take $f \in L^1 \cap \go$ and $\eps > 0$. By (\ref{main-20}),
$\exists A \in \sca_f$ with $\mu(A) \ge \| f \|_1 / \eps$. Set $g_\eps
= 1_A / \mu(A)$. We have that
\begin{equation}
  \label{main-30}
  \left| \frac{ \mu( (f \circ T^n) g_\eps) } {\mu(g_\eps)} \right|  = 
  | \mu( (f \circ T^n) g_\eps) | \le \| f \|_1 \, \| g_\eps \|_\infty 
  \le \eps.
\end{equation}
By Lemma \ref{lem-main},
\begin{equation}
  \limsup_{n \to \infty} \left| \frac{ \mu( (f \circ T^n) g) } {\mu(g)}
  \right| \le \eps
\end{equation}
holds for \emph{all} $g \in \lo$ with $\mu(g) \ne 0$.  Since $\eps$ is
arbitrary, we get that the above r.h.s.\ is zero. The case $\mu(g) =
0$ is trivial because the same assertion comes directly from
{\bfseries (GLM1)}. This proves {\bfseries (LLM)}, namely, assertion
\emph{(b)}.

Finally for \emph{(c)}. Take a $G \in \go$ such that $\avg(G) >
0$. Since $\mu_V(G) \to \avg(G)$, as $\iv$, {\bfseries (GGM2)} implies
that there exist a large enough $M$ and a $V \in \scv$, with $\mu(V)
\ge M$, such that
\begin{equation}
  \label{main-40}
  | \mu_V ((F \circ T^n)G) - \avg(F) \avg(G) | \le \eps \mu_V(G) 
\end{equation}
for all $n \ge M$. Setting $g := G 1_V$, we can divide (\ref{main-40})
by $\mu_V(G) = \mu(g) / \mu(V)$ and take the $\limsup$ in $n$:
\begin{equation}
  \limsup_{n \to \infty} \left| \frac{ \mu( (F \circ T^n) g) } {\mu(g)}
  - \frac{\avg(G)}{\mu_V(G)} \, \avg(F) \right| \le \eps.
\end{equation}
By Lemma \ref{lem-main}, the above holds $\forall g \in \lo$, with
$\mu(g) \ne 0$. Since $\eps$ can be taken arbitrarily close to 0 and
$\avg(G) / \mu_V(G)$ arbitrarily close to 1, we have that, for all $F
\in \go$ and $g \in \lo$, with $\mu(g) \ne 0$,
\begin{equation}
  \lim_{n \to \infty} \mu( (F \circ T^n) g)  = \avg(F) \mu(g).
\end{equation}
The corresponding statement for $\mu(g) = 0$ comes from
{\bfseries (GLM1)}.
\qed

\section{The equilibrium observables}
\label{sec-eo}

The ``pure'' \erg\ theorist might raise an eyebrow at the
constructions of Section \ref{sec-setup}, especially at the ideas of
the exhaustive family (which demands that one singles out some sets as
more important than the others) and of the infinite-volume average
(which is not a \me, or even guaranteed to always exist).

Though these issues (and more) have been addressed in \cite{limix},
one might still want to see if some of the concepts presented here can
be viewed from the vantage point of traditional infinite \erg\ theory.
For what follows I am indebted to R.~Zweim\"uller.

As we discussed in the introduction, the definition {\bfseries (GLM2)}
makes sense as a kind of convergence to equilibrium for a large class
of initial distributions (see also the observation on {\bfseries
(GLM1)} after the statement of Theorem \ref{thm-main}). Without
worrying too much about predetermining good test \fn s for this
convergence (namely, the global \ob s), and the value of any such
limit (namely, the infinite-volume average), one might simply consider
the space $\eo = \eo(\ps, \sca, \mu, T)$ of \emph{all} the good test
\fn s, in this sense:
\begin{equation}
  \label{def-eo}
  \eo := \rset{F \in L^\infty} {\exists \rho(F) \in \R \
  \mbox{s.t.} \ds \lim_{n \to \infty} \mu((F \circ T^n) g) = 
  \rho(F) \mu(g), \, \forall g \in L^1 }.
\end{equation}
(Occasionally, one might want to restrict the space of the initial
distributions to some subspace of $L^1$.) Clearly, $\eo$ is a vector
space which contains at least the constant \fn s.

$\rho(F)$ represents a sort of \emph{value at equilibrium} of $F$ and,
in this context, it need not have anything to do with $\avg(F)$ (which
might or might not exist), $\scv$, or the choice of $\go$ and $\lo$.
Thus, the elements of the vector space $\eo$ may be called the
{\bfseries equilibrium \ob s} and $\rho: \eo \into \R$ the {\bfseries
equilibrium functional}.

If we are in either case (H1) or (H2), Theorem
\ref{thm-main}\emph{(d)} shows that, for a given $F \in \go$, one only
need find \emph{one} local \ob\ that verifies the limit in
(\ref{def-eo}). Also, by Theorem \ref{thm-main}\emph{(b)}, any $f \in
\go \cap L^1$ belongs to $\eo$, with $\rho(f) = 0$. Therefore, in
these cases, it makes sense to introduce $\hat{\eo} := \eo / (\go \cap
L^1)$, and $\rho$ is well defined there.  When talking about
$\hat{\eo}$, we write $F \in \hat{\eo}$ to mean $F \in \eo$, and $F =
G$ to mean $[F] = [G]$ (where $[ \cdot ]$ denotes an equivalence class
in $\eo / (\go \cap L^1)$).

Determining $\hat{\eo}$ for a given, say, exact \dsy\ appears to be as
complicated as proving {\bfseries (GLM2)} for a truly large class of
global \ob s, though occasionally some information can be obtained
quickly. We conclude this note by giving some examples thereof.

\bigskip
\noindent
\textit{Boole transformation.} This is the transformation $T: \R \into
\R$ defined by $T(x) := x - 1/x$. This map preserves the Lebesgue \me\
on $\R$, as it is easy to verify, and is exact \cite{a}. We can use
the fact that $T$ is odd to construct a nonconstant equilibrium
\ob. Set $F(x) := \mathrm{sign}(x)$, and $g := 1_{[-1,1]}$. Clearly,
for all $n \in \N$, $F \circ T^n$ is odd and $\mu((F \circ T^n)g) =
0$, so $F \in \hat{\eo}$, with $F \ne $ constant, and $\rho(F) = 0$.

Evidently, the same reasoning can be applied to any exact map with an
odd symmetry.

\bigskip
\noindent
\textit{Translation-invariant expanding maps of $\R$}.  Take a $C^2$
bijection $\Phi: [0,1] \into [k_1, k_2]$, with $k_1, k_2 \in \Z$, and
$\Phi' > 1$, where $\Phi'$ denotes the derivative of $\Phi$. (Notice
that these conditions imply $\Phi(0) = k_1$, $\Phi(1) = k_2$, and $k
:= k_2- k_1 \ge 2$.) Define $T: \R \into \R$ via
\begin{equation}
  \label{eo-10}
  T |_{[j, j+1)} (x) := \Phi(x - j) + j,
\end{equation}
for all $j \in \Z$. By construction $T(x + 1) = T(x) +1$, $\forall x
\in \R$, and so $T$ is a $k$-to-1 translation-invariant map, in the
sense that it commutes with the natural action of $\Z$ in $\R$.

Suppose that $T$ preserves the Lebesgue \me, which we denote $m_\R$.
(One can easily construct a large class of maps of this kind.)  It can
be proved that any such $T$ is exact \cite{lmixmap}.  Now, define $I
:= [0,1)$ and $T_I : I \into I$ as $T_I (x) := T(x)$ mod 1.  Clearly,
$(I, \scb_I, T_I, m_I)$, where $\scb_I$ and $m_I$ are, respectively,
the Borel $\sigma$-algebra and the Lebesgue \me\ on $I$, is a
\pr-preserving \dsy. It is easy to see that it is exact, and thus
mixing.

Now consider a $\Z$-periodic, bounded, $F: \R \into \R$. Evidently,
$\forall x \in I$, $\forall n \in \N$, $F \circ T^n (x) = F \circ
T_I^n (x)$. Hence, by the mixing of the quotient \dsy, for any
square-integrable $g$ supported in $I$,
\begin{equation}
\label{eo-20}
\begin{split}
  \lim_{n \to \infty} m_\R ((F \circ T^n) g) &= \lim_{n \to \infty}  
  m_I ((F \circ T_I^n) g) \\ 
  &= m_I(F) \, m_I(g)  \\
  &= m_I(F) \, m_\R(g).
\end{split}
\end{equation}
By the exactness of $T$, the above holds for all $g \in
L^1(\R)$. Hence $F \in \hat{\eo}$, with $\rho(F) = m_I(F) =
\overline{m_\R}(F)$.

An analogous procedure (using $I_j := [0,j)$ instead of $I$) can be
employed to prove that any $(j\Z)$-periodic, bounded $F$ belongs in
$\hat{\eo}$, with $\rho(F) = \overline{m_\R} (F)$.  In
\cite{lmixmap} we extend this result to \ob s that are
quasi-periodic w.r.t.\ any $j\Z$, and more.

\bigskip
\noindent
\textit{Random walks}. A special case of the above situation occurs
when $\Phi$ is linear. The result is a piecewise linear Markov map
that represents a random walk in $\Z$, in the following sense. Denote
by $\lfloor x \rfloor$ the maximum integer not exceeding $x \in
\R$. If an initial condition $x \in I$ is randomly chosen with law
$m_I$, then the stochastic process $( \lfloor T^n(x) \rfloor )_{n \in
\N}$ is precisely the random walk starting in $0 \in \Z$, with
uniform transition probabilities for jumps of $k_1, k_1 +1, \ldots,
k_2 - 1$ units \cite{liutam}.

A reelaboration of a result of \cite{limix} shows that $\hat{\eo}$
contains all $L^\infty$ \fn s such that the limit
\begin{equation}
  \label{eo-30}
  \rho (F) := \lim_{M \to \infty} \int_{a-M}^{a+M} \!\! F(x) \, dx
\end{equation}
exists independently of and uniformly in $a \in \R$.  In fact, it is
proved in \cite[Thm.~4.6(b)]{limix} (see also \cite[Thm.~9]{liutam})
that, 
if $g \in L^1$, 
$F \in L^\infty( \sca_0 )$, where $\sca_0$ is the
$\sigma$-algebra generated by the partition $\{ [j, j+1) \}_j$,
and the limit
\begin{equation}
  \label{eo-40}
  \lim_{j \to \infty} \int_{q-j}^{q+j} \!\! 
  F(x) \, dx =: \overline{m_\R}(F)
\end{equation}
($j \in \Z$) exists uniformly in $q \in \Z$, then $m_\R ((F \circ
T^n)g) \to \overline{m_\R}(F) m_\R(g)$, as $n \to \infty$. Obviously,
comparing (\ref{eo-30}) with (\ref{eo-40}), $\rho(F) =
\overline{m_\R}(F)$.

Now, for a general $F$, one can take $g = 1_{[0,1)} \in
L^1(\sca_0)$. It is easy to check that $P^n g$ is $\sca_0$-measurable
too, thus
\begin{equation}
\label{eo-50}
\begin{split}
  \lim_{n \to \infty} m_\R ((F \circ T^n) g) &= \lim_{n \to \infty}  
  \la \E(F | \sca_0), P^n g \ra \\
  & = \overline{m_\R}(\E(F | \sca_0)) \, m_\R(g) \\
  &=  \rho(F) \, m_\R(g),
\end{split}
\end{equation}
which proves our claim.

If the random walk has a drift, say a positive drift, then
a.e.\ \o\ will converge to $+\infty$. Therefore, any bounded \fn\ $G$
that asymptotically shadows any of the above \ob s---meaning $\lim_{x
\to +\infty} (G(x) - F(x)) = 0$, for some $F$ verifying
(\ref{eo-30})---will also belong to $\hat{\eo}$, with $\rho(G) =
\rho(F)$.

\end{document}